\documentclass[amssymb]{article}
\usepackage{amstex}
\usepackage{amssymb}

\newcounter{letters}

\def\E{\ifmmode{\mathbb E}\else{$\mathbb E$}\fi} 
\def\N{\ifmmode{\mathbb N}\else{$\mathbb N$}\fi} 
\def\R{\ifmmode{\mathbb R}\else{$\mathbb R$}\fi} 
\def\Q{\ifmmode{\mathbb Q}\else{$\mathbb Q$}\fi} 
\def\C{\ifmmode{\mathbb C}\else{$\mathbb C$}\fi} 
\def\H{\ifmmode{\mathbb H}\else{$\mathbb H$}\fi} 
\def\Z{\ifmmode{\mathbb Z}\else{$\mathbb Z$}\fi} 
\def\P{\ifmmode{\mathbb P}\else{$\mathbb P$}\fi} 
\def\T{\ifmmode{\mathbb T}\else{$\mathbb T$}\fi} 
\def\SS{\ifmmode{\mathbb S}\else{$\mathbb S$}\fi} 
\def\DD{\ifmmode{\mathbb D}\else{$\mathbb D$}\fi} 

\title{Lectures on Seiberg-Witten Invariants}
\author{S.Akbulut}
\date{}
\begin{document}
\maketitle
\def \Di {D\!\!\!\!/\,}
\def \partiali {\partial\!\!\!/\,}

In October 1994 Seiberg-Witten invariants entered in $4$-manifold theory
with a big
bang. Not only  did these invariants tidy up the  Gauge Theory, but they
also gave
some exciting new results on topology of smooth $4$-manifolds. These notes
grew out
of the lectures I have given in learning seminars at MPI in Bonn, and METU
in Ankara
on this subject. The main goal of these notes is not to survey the whole
area, but
rather establish conventions for novice topologist like myself, and go
through some
recent selected results. In these notes I avoided the general Clifford algebra
constructions in favor of more direct representation theory of
$Spin_{c}(4)$.

 I have benefited greatly from stimulating papers [{\bf KM}], [{\bf W}],
and [{\bf
T}], as well as unpublished lecture notes of C.Taubes. I also benefited
seminar
talks  by  R.Fintushel, T.Parker, and T.Draghici at MSU. I thank I.Hambelton
and
T.Onder for inviting me to MPI and METU, giving me the chance to work on these
notes. I make no claim of originality in these  notes, they are merely  modest
efforts to understand the  results  from the sources mentioned above.

\section{Introduction}

Every compact oriented smooth $4$-manifold has a $Spin_{c}$ structure, i.e. the
second Steifel-Whitney  $w_{2}(X)\in H^{2}(X;\Z_{2})\;$ has an integral
lifting. This is because:
$w_{2}(X)$ can be represented by an imbedded surface $F\subset X$. If $F$ is
orientable then clearly the homology class
$[F]$ comes from an integral class; if not then it suffices to show the circle
$S\subset F$ representing $w_{1}(F)$ is null homologous in $H_{1}(X;\Z)$,
because the Bockstein $\delta[F]=[S]$ in the coefficient exact sequence:
$$ ..\to H_{2}(X;\Z)\stackrel{\times 2}{\longrightarrow}
H_{2}(X;\Z)\stackrel{\rho}{\longrightarrow}
H_{2}(X;\Z_{2})\stackrel{\delta}{\to}
H_{1}(X;\Z)\to..$$ where $\rho $ is the reduction map. Now if
$\delta[F]\neq 0$, we
can choose an imbedded oriented $3$-manifold $\Sigma \subset X$ representing
the
Poincare dual of
$\delta[F]$,
 which is transverse to $F$. Then  $T=F\cap \Sigma \subset F$ has a trivial
normal
bundle $\nu (T,F)$ since
$$\nu (T,X) = \nu (T,F)\oplus \nu (T,\Sigma)$$ and the two other normal
bundles in
the above equality are trivial. This gives a contradiction, since in $F$ the
$1$-manifold $T$ meets $S$ transversally at one point and $[S]=w_{1}(F)$
implies
$\nu (T,F)$ must necessarily be nontrivial $\;\;\;\Box$.
\vspace{.05in}
\begin{eqnarray*}\mbox{Recal:}\;\;\;\;\;\;\; Spin(4)&=&SU (2)\times SU (2)\\
Spin_{c}(4)&=&(\;SU (2)\times SU(2)\times S^{1}\;)/\Z_{2} =(\;Spin(4)\times
S^{1})/\Z_{2}\\ SO(4)&=&(\;SU (2)\times SU (2)\;)/\Z_{2}\\ U(2)&=&(\;SU
(2)\times S^{1}\;)/\Z_{2}\end{eqnarray*}

\noindent We have fibrations:
$$S^{1}\longrightarrow Spin_{c}(4)\to SO(4)$$
$$\Z_{2}\longrightarrow Spin_{c}(4)\to SO(4)\times S^{1}$$

\vspace{.01in}

\noindent We can also identify
$Spin_{c}(4)=\{ (A,B)\in U_{2}\times U_{2}\; |\; det(A)=det(B)\;\}$ by

 $$(A,B)\leadsto (A.(det A)^{-1/2}\; I\;,\; B.(det B)^{-1/2}\; I\;,\; (det
A)^{1/2}
)$$

We also have $2$ fold cover $Spin_{c}(4)\to SO(4)\times S^{1}$ .The
fibrations above
extend to fibrations:
$$S^{1}\to Spin_{c}(4)\to SO(4)\to K(\Z,2)\to BSpin_{c}(4)\to BSO(4)\to
K(\Z,3)$$ The last map in the sequence is given by the Bokstein of the second
Steifel-Whitney class
$\delta (w_{2})$ which explains why lifting of $w_{2}$ to an integral class
corresponds to a $Spin_{c}(4)$-structure. We also have the fibration:
$$\Z_{2}\to Spin_{c}(4)\to SO(4)\times S^{1}\to  K(\Z_{2},1)\to
BSpin_{c}(4)\to BSO(4)\times BS^{1}\to K(\Z_{2},2) $$ The last map in this
sequence is given by
$\;w_{2}\times 1 + 1\times
\rho(c_{1})\;$ which clearly vanishes exactly when  $\delta (w_{2})=0$ .
Finally we
have the fibration:
$$\Z_{2}\to Spin(4)\times S^{1}\to Spin_{c}(4) \to  K(\Z_{2},1)\to
BSpin(4)\times BS^{1}\to BSpin_{c}(4)\to K(\Z_{2},2) $$

\noindent The last map is given by $w_{2}$. This sequence says that locally a
$\;Spin_{c}(4)$ bundle consists a pair of a $Spin(4)$ bundle and a complex line
bundle. Also recall
$\;\;H^{2}(X;\Z)=[X,K(\Z,2)]=[X,BS^{1}]= \{\mbox{complex line bundles on
X}\}$

\vspace{.15in}

\noindent {\bf Definition}: Let $ L\longrightarrow X$ be a complex line
bundle  over
a smooth oriented $4$-manifold with
$c_{1}(L)=w_{2}(TX)$ (i.e. $L$ is a characteristic line bundle). A
$Spin_{c}(4)\;$ structure on $X$, corresponding $L$, is a principal
$\;Spin_{c}(4)$-bundle
$\;P\longrightarrow X\;$ such that the associated framed bundles of
$TX$ and $L$  satisfy:
$$P_{SO(4)}(TX)=P\times_{{\rho}_{0}} SO(4) $$
$$P_{S^{1}}(L)=P\times_{{\rho}_{1}} S^{1} $$

\noindent where $(\rho_{0}, \rho_{1}): Spin_{c}(4)\to SO(4)\times S^{1}$ are
the
obvious projections

$$\begin{array}{ccccc} &  & Spin(4)\times S^{1} &  & \\ &&&&\\ &  &
\downarrow \pi & &\\ &&&&\\ SO(4) & \stackrel{\rho_{0}} {\longleftarrow }
&Spin_{c}(4)&
\stackrel{\rho_{1}} { \longrightarrow} & S^{1}\\ &&&&\\ &&&&\\
 & \swarrow \rho_{+} && \rho_{-}\searrow & \\ &&&&\\
 U(2)  && \downarrow
\tilde{\rho}_{+}\;\;\;\;\;\;\tilde{\rho}_{-}\downarrow & & U(2)\\ &&&&\\
 & Ad \searrow&&  \swarrow Ad& \\ &&&&\\ &  &SO(3) & & \\

\end{array}$$

So $\;\tilde{\rho}_{\pm}= Ad\circ \rho _{\pm}\;$, also call
$\;\bar{\rho}_{\pm}=
\rho_{\pm}\circ \pi\;$.
 For $x\in \H = \R^{4}$ we have
\begin{eqnarray*}
\rho_{1}[\;q_{+},q_{-},\lambda \;]&=&\lambda^{2} \\
\rho_{0}[\;q_{+},q_{-},\lambda
\;]&=&[\;q_{+},q_{-}\;]\;\;\;\;\;,\;\mbox{i.e.}
\;\;\;x\longmapsto q_{+}xq_{-}^{-1}\\
\rho_{\pm}[\;q_{+},q_{-},\lambda \;]&=&[\;q_{\pm},\lambda \;]
\;\;\;\;\;\;\;,\;\mbox{i.e.}
\;\;\;x\longmapsto q_{\pm} x \lambda^{-1} \\
\tilde{\rho}_{\pm}[\;q_{+},q_{-},\lambda \;]&=& Ad\circ q_{\pm}
\;\;\;\;\;\;\;,\;
\mbox{i.e.}\;\;\;x\longmapsto q_{\pm}xq_{\pm}^{-1}\\
\bar{\rho}_{\pm}(\;q_{+},q_{-},\lambda \;) &=& \lambda q_{\pm}
\end{eqnarray*}

\noindent Apart from $TX$ and $L$, $Spin_{c}(4)$ bundle  $P\to X$ induces a
pair of
$U(2)$ bundles:
$$W^{\pm}=P\times _{\rho_{\pm}}\C^{2}\longrightarrow X$$ Let
$\Lambda^{p}(X)=\Lambda^{p}T^{*}(X) $ be the bundle of exterior $p$ forms.
If $X$ is
a Riemanian manifold (i.e. with metric), we can construct the bundle of
self(antiself)-dual 2-forms
$\;\Lambda_{\pm}^{2}(X)$ which we abbreviate by $ \Lambda^{\pm}(X)
\;$.  We can identify   $ \Lambda ^{2}(X) \;$ by the Lie algebra
$so(4)$-bundle
$$\Lambda ^{2}(X)=P(T^{*}X)\times _{ad}so(4)\;\;\;\;\;\;\mbox{by}\;\;\;
\Sigma \;a_{ij}\;dx^{i}\wedge dx^{j}\;\longleftrightarrow \; (a_{ij})$$ where
$ad:SO(4)\to so(4)\;$ is the adjoint representation. The adjoint action
preserves
the two summands of
$so(4)=spin(4)=so(3)\times so(3)=\R^{3}\oplus \R^{3}$. By above
identification it is easy to see that the $\pm 1$ eigenspaces $
\Lambda^{\pm}(X) \;$
of the star operator
$*:\Lambda (X)\to \Lambda (X)$ corresponds to these two ${\R}^{3}$-bundles;
this
gives:
$$\Lambda^{\pm}(X)= P\times _{\tilde{\rho}_{\pm}} \R^{3}  $$

 If the $Spin_{c}(4)$ bundle $ P\to X$ lifts to $Spin(4)$ bundle $
\bar{P}\to X$ (i.e. when $w_{2}(X)=0$), corresponding to the obvious
projections $
p_{\pm}: Spin(4)\to SU(2) $,
$p_{\pm}(q_{-},q_{+})=q_{\pm}$ we get a pair of $SU(2)$ bundles:
$$V^{\pm}= P\times_{p_{\pm}} \C^{2}  $$ Clearly since
$\;x\longmapsto q_{\pm} x  \lambda ^{-1}=
 q_{\pm} x \;(\lambda^{2})^{-1/2}\; $ in this case we have:
$$\;W^{\pm}=  V^{\pm}\otimes  L^{-1/2}\;$$

\subsection {Action of $\;\Lambda^*(X) \;\;$on$\;\; W_\pm$ }

 From the definition of $Spin_{c}(4)$ structure above we see that
$$T^{*}(X)=P\times \H/ (p,v)\sim
(\tilde{p},q_{+}v\;q_{-}^{-1})\;\;,\;\mbox{where}\;\;\;
\tilde{p}=p[\;q_{+},q_{-},\lambda\;]$$ We define left actions (Clifford
multiplications), which is  well defined by
$$  T^{*}(X) \otimes  W^{+}\longrightarrow W^{-}\;\;,\;\mbox{by}\;\;\;\;
[\;p,v\;]\otimes[\;p,x\;]\longmapsto [\;p,-\bar{v}x\;]$$
$$T^{*}(X) \otimes  W^{-} \longrightarrow W^{+}\;\;,\;\mbox{by}\;\;\;\;
[\;p,v\;]\otimes[\;p,x\;]
\longmapsto [\;p,vx\;]$$  From identifications, we can check the well
definededness
of these actions, e.g.:
 $$[p,v\;]\otimes [\;p,x\;]\sim [\;\tilde{p},q_{+}v\;q_{-}^{-1}\;]
\otimes [\;\tilde{p}, q_{+}x\lambda^{-1}\;]
\longmapsto [\;\tilde{p}, \; q_{-}(-\bar{v}x)\lambda^{-1}\;]
\sim[\;p,-\bar{v}x\;]$$
\vspace{.05in}

\noindent By dimension reason complexification of these representation give
$$ \rho:T^{*}(X)_\C\stackrel{\cong }{\longrightarrow} Hom(W^{\pm},
W^{\mp})\equiv W^{\pm}\otimes W^{\mp}$$

\noindent We can put them together as a single representation (which we
still call
$\rho $)
$$ \rho: T^{*}(X)\longrightarrow Hom(W^{+}\oplus W^{-})\;\;,\;\;\mbox{by}
\;\;\;v\;\longmapsto \;
\rho (v)=\left(\begin{array}{cc} 0& v \\
 -\bar{v} & 0\end{array} \right) $$ We have  $\;\rho (v)\circ\rho
(v)=-|v|^{2}I\;$.
By  universal property of the Clifford algebra this representation extends
to the
Clifford algebra
$C(X)=\Lambda^{*}(X)$ (exterior algebra)

$$\begin{array}{ccc}
\Lambda^{*}(X)& & \\ & & \\
\downarrow &\searrow & \\ & & \\ T^{*}(X) & \longrightarrow & Hom\;(W^{+}\oplus
W^{-})
\end{array}$$

\noindent One can  construct this extension without the aid of the universal
property of the Clifford algebra, for example since
$$\Lambda^{2}(X)=\left \{\; v_{1}\wedge v_{2}=\frac{1}{2}(v_{1}\otimes
v_{2}-v_{2}\otimes v_{1})\;|\; v_1,v_{2}\in T^{*}(X)\;\right \}$$ The action of
$T^{*}(X)$ on
$W^{\pm}$ determines the action of
$\;\Lambda^{2}(X)=\Lambda^{+}(X)\otimes\Lambda^{-}(X) $, and since
$\;2Im\;(v_{2}\bar{v}_{1})=-v_{1}\bar{v}_{2} + v_{2}\bar{v}_{1} $ we have
the action
$\rho$ with property:
$$\; \Lambda^{+}(X)\otimes  W^{+} \longrightarrow W^{+}\;\;\;\;\mbox{to
be}\;\;\;
[\;p,v_{1}\wedge v_{2}\;]\otimes [\;p,x\;]
\longrightarrow [\;p,Im\;(v_{2}\bar{v}_{1})x\;] $$
 $$\rho: \Lambda^{+} \longrightarrow Hom(W^{+}, W^{+})$$
\begin{eqnarray}\rho (v_{1}\wedge v_{2})&=&
\;\frac{1}{2}\;[\;\rho(v_{1}),\rho (v_{2})\;]
\end{eqnarray}
\vspace{.1in} Let us write the local descriptions of these representations:
We first
pick a
 local orthonormal basis $\;\{e^{1},e^{2},e^{3},e^{4}\}\;$  for
$T^{*}(X)$, then we can take
$$\;\{\;f_{1}=\frac{1}{2}( e^{1}\wedge e^{2} \pm e^{3}\wedge e^{4}),\;
f_{2}=\frac{1}{2}( e^{1}\wedge e^{3} \pm e^{4}\wedge e^{2}),\;
f_{3}=\frac{1}{2}(
e^{1}\wedge e^{4} \pm e^{2}\wedge e^{3})\;\}$$ to be  a basis for
$\Lambda^{\pm}(X)$. After the local identification $T^{*}(X)=\H$ we can take
$e^{1}=1,\;e^{2}=i,\;e^{3}=j,\; e^{4}=k$. Let us identify
$W^{\pm}=\C^{2}=\{z+jw\;|\; z,w\in \C\;\}$, then the multiplication by
$1,i,j,k$ (action on $\C^{2}$ as multiplication on left) induce the
representations
$\rho (e^{i})\;,\;i=1,2,3,4$.  From this we see that
$\Lambda^{+}(X)$ acts trivially on $W^{-}$; and the basis
$f_{1},f_{2},f_{3}$ of
$\Lambda^{+}(X)$ acts on $W^{+}$ as multiplication by $i,j,k$, respectively
(these
are called Pauli matrices).

$$\begin{array}{cc}
\rho(e^{1})=\left( \begin{array} {cccc}  &  &1 & 0\\  &  & 0 & 1
\\-1 & 0 &  &  \\ 0 & -1 &  &  \end{array}
\right) & \;\;\;\;
\rho(e^{2})=\left( \begin{array}{cccc} &&i&0\\ &&0&-i \\i&0&&\\ 0&-i&&
\end{array} \right) \\ &\\ &\\
\rho(e^{3})=\left( \begin{array}{cccc} &&0&-1\\ &&1&0 \\0&-1&&\\1&0&&
\end{array} \right) &\;\;\;
 \rho(e^{4})=\left( \begin{array}{cccc} &&0&-i\\ &&-i&0
\\0&-i&&\\-i&0&&
\end{array} \right) \\
\end{array}$$

\vspace{.15in}

$$\begin{array}{ccc}
\rho(f_1)=\left( \begin{array}{cc} i&0\\ 0&-i \end{array} \right)&
\rho(f_2)=\left( \begin{array}{cc} 0&-1\\ 1&0 \end{array} \right)&
\rho(f_3)=\left( \begin{array}{cc} 0&-i\\ -i&0 \end{array} \right)
\\
\end{array}$$

\vspace{.15in}

In particular we get an isomorphism $\Lambda^{+}(X)\longrightarrow su\;(W^{+})$
(traceless skew adjoint endemorphism of $W^{+}$); which after complexifying
extends
to an isomorphism
$\rho:\Lambda^{+}(X)_\C\cong sl\;(W^{+})$ (traceless endemorphism of $W^{+}$)

$$\begin{array}{ccc}
\Lambda^{+}(X)&\stackrel{\cong}{\longrightarrow}& su\;(W^{+})\\ &&\\
\bigcap & & \bigcap\\ &&\\
\Lambda^{+}(X)_\C&\stackrel{\rho}{\longrightarrow} & sl\;(W^{+})

\end{array}$$

\vspace{.1in}

Recall $\;Hom(W^{+}, W^{+})\cong W^{+}\otimes (W^{+})^{*}\;$; we identify
the dual
 space $(W^{+})^{*}$ naturally with $\bar{W}^{+}$ (= $W^{+}$ with  scalar
multiplication  $\;c.v=\bar {c}v$) by the pairing
$$W^{+}\otimes \bar{W}{+}\longrightarrow \C $$ given by
$\;z\otimes w\to z\bar{w}$. Usually
$sl\;(W^{+})$ is  denoted by $(W^{+}\otimes \bar{W}^{+})_0\;$ and the trace map
gives the identification:
$$W^{+}\otimes \bar{W}^{+}=(W^{+}\otimes \bar{W}^{+})_0\oplus \C=
\Lambda^{+}(X)_\C\oplus \C$$ Let  $\sigma
:W^{+}\longrightarrow\Lambda^{+}(X)\; $ be the map
$\;[\;p,x\;]\longmapsto [\;p,\;\frac{1}{2}(x i \;\bar{x})\;] $. By local
identification as above $ W^{+}=\C^{2}$ and
$\Lambda^{+}(X)=\R\oplus \C$,  we see  $\sigma$ corresponds to

$$\;(z,w)\longmapsto i\;\left(\frac{|z|^{2}-|w|^{2}}{2}\right) -k\;
Re(z\bar{w})+
j\;Im(z\bar{w}) =
\left(\frac{|w|^{2}-|z|^{2}}{2}\right) + z\bar{w} $$

\noindent  We identify this by the element $i\sigma(z,w)$ of
$su\;(W^{+})$ (by  Pauli matrices) where:
\begin{eqnarray}\;(z,w)\; \longmapsto \sigma(z,w)\;=\left(
\begin{array}{cc} (\;|z|^{2}-|w|^{2})/2& {z}\bar{w}\\
\bar{z}{w}& (\;|w|^{2}-|z|^{2})/2 \end{array} \right)\end{eqnarray}
 $\sigma$ is the projection of the diagonal elements of $W^{+}\otimes
\bar{W}^{+}$ onto
$(W^{+}\otimes \bar{W}^{+})_0$

\vspace{.05in}

\noindent We can check:
\begin{eqnarray} i\;\sigma(z,w)=\rho
\;[\;\frac{|z|^{2}-|w|^{2}}{2}\; f_{1} + Im(z\bar{w})
\;f_{2} - Re(z\bar{w})\; f_{3}\;]
\end{eqnarray}

 From these identifications we see:
\begin{eqnarray} |\;\sigma(\psi) \;|^{2}&=&\frac{1}{4} |\;\psi
\;|^{4} \\ <\;\sigma(\psi)\;\psi,\psi\;> &=& \frac{1}{2} |\;\psi
\;|^{4}\\
 <\rho(\omega)\;\psi\;,\;\psi>&=&2i\;<\rho (\omega)\;,i\;
\sigma(\psi)>
\end{eqnarray} Here the norm in $\;su(2)\;$ is induced by the inner product
$<A,B>=\frac{1}{2}trace(AB)\;$. \\

 By calling $\sigma(\psi, \psi)=\sigma (\psi)$ we  extend the definition of
$\sigma $ to $ W^{+}\otimes \bar{W}^{+} $ by
$$<\rho(\omega)\;\psi\;,\;\varphi>=2i\;<\rho (\omega)\;,i\;
\sigma(\psi,
\varphi)> $$
\vspace{.02in}
$$\begin{array}{cccc}
\Lambda^{+}(X)\;\;\;\;\;=&su (W^{+})&\;\;\stackrel{i\;\sigma}{\longleftarrow}&
W^{+}\\ &&&\\ &\bigcap & & \bigcap\\ &&&\\ (W^{+}\otimes \bar{W}^{+})_0\;\;
=&sl\;(W^{+})&\stackrel{i\;\sigma}{\longleftarrow} & W^{+}\otimes
\bar{W}^{+}
\end{array}$$

\vspace{.15in}

\noindent{\bf Remark}: A $Spin_{c}(4)$ structure can also be defined as a
pair of
$U(2)$ bundles:
$$W^{\pm}\longrightarrow
X\;\;\mbox{with}\;\;det(W^{+})=det(W^{-})\longrightarrow X
\;\;\mbox{(a complex line bundle), }$$
$$\mbox{ and an action}\;\; c_{\pm}:T^{*}(X)\longrightarrow
Hom(W^{\pm},W^{\mp})\;\;
\mbox{with}\;
\;c_{\pm}(v)c_{\mp}(v)=-|v|^{2}I $$

\vspace{.1in}

\noindent The first definition clearly implies this, and conversely we can
obtain
the first definition by letting the principal
$Spin_{c}(4)$ bundle to be:
$$P=\{\;(p_{+},p_{-})\in P(W^{+})\times P(W^{-})\;|\;
det(p_{+})=det(p_{-})\;\} $$
Clearly, $\;Spin_{c}(4)=\{ (A,B)\in U_{2}\times U_{2}\; |\;
det(A)=det(B)\;\}$ acts
on
$P$ freely.

\vspace{.15in}

\noindent This definition  generalizes and gives way to the following
definition:

\vspace{.1in}

\noindent{\bf Definition}: A Dirac bundle  $W\longrightarrow X$ is a Riemanian
vector bundle with an action $
\rho:T^{*}(X)\longrightarrow Hom(W,W)\;$ satisfying
$\;\rho (v)\circ \rho (v)=-|v|^{2}I $.  $W$ is also equipped with  a
connection $D$
satisfying:
 $$<\rho (v)x,\rho (v)y>=<x,y>$$
  $$D_{Y}(\rho (v)s)=\rho (\nabla_{X} v) s + \rho (v) D_{Y}(s)$$
\noindent where $\nabla $ is the Levi-Civita connection on
$T^{*}(X)$, and $Y$ is a vector filed on $X$

\vspace{.1in}

An example of a Dirac bundle is $\;W=W^{+}\oplus W^{-}\longrightarrow X\;$ and
$\;D=d+d^{*}\;$ with
$W^{+}=\oplus\Lambda^{2k}(X)\; $ and
$W^{-}=\oplus\Lambda^{2k+1}(X)\; $ where
$\rho (v)= v\wedge . + v \;\bot ..$ (exterior $+$ interior product with
$v$). In this case
 $\rho : W^{\pm}\to W^{\mp} $. In the next section we will discuss the
natural
connections
$D$ for $Spin_{c}$ structures $W^{\pm}$ .

\section{ Dirac Operator}

Let $\cal{A}(L)$ denote the space of connections on a $\;U(1)$ bundle
$L\longrightarrow X$. Any $A\in \cal{A}(L)$ and the Levi-Civita connection
$A_{0}$
on the tangent bundle coming from Riemanian metric of
$X$ defines a product connection on $P_{SO(4)}\times P_{S^{1}}$. Since
$Spin_{c}(4)$ is the two fold covering of
$SO(4)\times S^{1}$, they have the same Lie algebras
$spin_{c}(4)=so(4)\oplus i\;\R$. Hence we get a connection
$\tilde{A}$ on the $Spin_{c}(4)$ principle bundle
$P\longrightarrow X$. In particular the connection
$\tilde{A}$ defines connections to all the associated bundles of P, giving back
$A,\;A_{0}$ on $L, \;T(X)$ respectively, and two new connections $A^{\pm}$ on
bundles $W^{\pm}$. We denote the corresponding covariant derivatives by
$\nabla_{A}$.
$$\nabla _{A} :\Gamma(W^{+})\to \Gamma(\;T^{*} X\otimes W^{+})$$ Composing
this with
the Clifford multiplication $\Gamma(\;T^{*} X\otimes W^{+} )\to \Gamma(W^{-}) $
gives the Dirac operator
$$\Di _{A} :\Gamma(W^{+})\to \Gamma(W^{-})$$ Locally, by choosing orthonormal
tangent vector field
$e=\{e_{i}\;\}_{i=1}^{4}$ and the dual basis of $1$-forms
$\{e^{i}\;\}_{i=1}^{4}$ in a neighborhood
$U$ of a point $x\in X$ we can write
$$\Di _{A}=\sum \rho (e^{i}) \nabla_{e_{i}}$$ where $\nabla_{e_{i}}:
\Gamma(W^{+})\to \Gamma(W^{+})$ is the covariant derivative
$\nabla _{A}$ along $e_{i}$. Also locally $W^{\pm}=V^{\pm}\otimes L^{1/2}
$, hence by Leibnitz rule, the connection  $A$ and the untwisted Dirac operator
$$\partiali : \Gamma(V^{+})\to \Gamma(V^{-})$$ determines $\Di_{A}$. Notice
that as
in $W^{\pm}$, forms $\Lambda^{*}(X)$ act on $V^{\pm}$. Now let
$\omega=(\omega_{ij})$ be the Levi-Civita connection $1$-form, i.e.
$so(4)$-valued ``equivariant"
$1$-form on $P_{SO(4\;)}(X)$ and
$\tilde{\omega } =(\tilde{\omega }_{ij})=e^{*}(\omega) $ be the pull-back
$1$-form on $U$. Since
$P_{SO(4\;)}(U)=P_{Spin(4\;)}(U)$ the orthonormal basis $e\in P_{SO(4)}(U)$
determines an orthonormal basis
$\sigma=\{\sigma^{k}\}\in P_{SU_{2}}(V^{+})$, then (e.g. [{\bf LM}])

$$\partiali(\sigma^{k})=\frac{1}{2}\sum_{i<j}
\rho( \tilde{\omega}_{ji})\;\rho( e^{i}) \rho ( e^{j})\;\sigma^{k} $$

Metrics on $T(X)$ and $L$ give metrics on $W^{\pm}$ and
$T^{*}(X) \otimes W^{\pm} $, hence we can define the adjoint $\nabla
_{A}^{*}:\Gamma(T^{*} X\otimes W^{-})\to
\Gamma(W^{+})$. Similarly we can define $\Di _{A} :\Gamma(W^{-})\to
\Gamma(W^{+}) $ which turns out to be the adjoint of the previous
$\Di_{A}$ and makes the following commute (vertical maps are Clifford
multiplications):
$$\begin{array}{ccccc}
\Gamma(W^{+})&\stackrel{ \nabla _{A}}{\longrightarrow}&\Gamma(\;T^{*} X
\otimes W^{+})&
\stackrel{ \nabla _{A}}{\longrightarrow}&\Gamma(\;T^{*} X \otimes T^{*}
X\otimes
W^{+})\\
\parallel & & \downarrow & &\downarrow \\
\Gamma(W^{+}) & \stackrel{ \Di _{A}}{\longrightarrow} &\Gamma(W^{-})&
\stackrel{ \Di _{A}}{\longrightarrow}&\Gamma(W^{+})\\
\end{array}$$

\noindent Let $F_{A}\in \Lambda^{2}(X)$ be the curvature of the connection
$A$ on  $L$, and
$F^{+}_{A}\in \Lambda^{+}(X)$ be the self dual part of this curvature, and
$s$ be
the scalar curvature of $X$. Weitzenbock formula says that:
\begin{eqnarray}
\Di_{A}^{2}(\psi)&=& \nabla _{A}^{*}\nabla _{A}\psi
 +\frac{s}{4}\psi + \frac{1}{4}\rho(F_{A}^{+})\psi
\end{eqnarray} To see this we we can assume
$\nabla_{e_{i}}(e^{j})=0$ at the point $x$
\begin{eqnarray*}
\Di_{A}^{2}\psi&=& \sum \rho( e^{i}).\;\nabla_{e_{i}}\;[\;\sum \rho
(e^{j}).\;\nabla_{e_{j}}\psi\;]\;\\ &=&
\nabla^{*}_{A}\nabla_{A}\psi+\frac{1}{2}\sum_{i,j}\rho (e^{i})\;\rho
(e^{j})\;(\nabla_{e_{i}}\nabla_{e_{j}} -
\nabla_{e_{j}}\nabla_{e_{i}})\;\psi\\ &=&
\nabla^{*}_{A}\nabla_{A}\psi+\frac{1}{2}
\sum_{i,j}\rho (e^{i})\;\rho (e^{j})\; \Omega_{ij}^{A} \;\psi
\end{eqnarray*}
$\Omega_{ij}^{A}=R_{ij}+ \frac{1}{2} F_{ij}\;$ is curvature on
$V^{+}\otimes L^{1/2}$, i.e. $R_{ij}$ is the Riemanian curvature and the
imaginary
valued $2$-form
$F_{ij}$ is the curvature of $A$ for the line bundle $L$ (endemorphisims of
$W^{+}$). So if $\psi=\sigma\otimes \alpha  \in
\Gamma(V^{+}\otimes L^{1/2})\;$, then

\begin{eqnarray*}\frac{1}{2}\sum_{i,j}\rho (e^{i})\;\rho (e^{j})\;
\Omega_{ij}^{A}\;(\sigma\otimes \alpha)&=&\frac{1}{2}
 (\sum \rho (e^{i})\;\rho (e^{j})\;R_{ij} \;\sigma )\otimes \alpha
\\ && + \;\frac{1}{4}\sum \rho (e^{i})\;\rho (e^{j}) \sigma \otimes ( F_{ij}
\alpha)\\ &&\\ &=& \frac{1}{8} \;\sum \rho (e^{i})\;\rho (e^{j})\;
\rho (e^{k})\;
\rho (e^{l})\;R_{ijkl} \;(\psi)  \\ && + \;\frac{1}{4}\;\rho
\;(\sum  F_{ij} \; e^{i}\wedge e^{j}) \;(\psi )
 \end{eqnarray*} The last identity follows from (1). It is a standard
calculation
that the first term is $s/4$ ( e.g.[{\bf LM}], pp. 161), and since $\Lambda
^{-}
(X)$ act as zero on $W^{+}$, the second term can be replaced by
$$\frac{1}{4}\;\rho(F_{A}^{+})\;\psi=\frac{1}{4}\;\rho \;(\sum  F_{ij}^{+}
\; e^{i}\wedge e^{j}) \;\psi $$

\subsection {A Special Calculation}

In Section 4 we need some a special case  (7). For this, suppose
$$V^{+}=L^{1/2}\oplus L^{-1/2}$$ where  $L^{1/2}\longrightarrow X$ is some
complex
line bundle with
$L^{1/2}\otimes L^{1/2}=L$.
 Hence $W^{+}=(L^{1/2}\oplus L^{-1/2})\otimes L^{-1/2}=L^{-1}\oplus \C$. In
this case there is a unique connection
$\frac{1}{2}A_{0}$ in $L^{-1/2}\to X$  such that the induced Dirac operator
$D_{A_{0}}$ on $W^{+}$ restricted to the trivial summand $\underline \C\to X$
is the exterior derivative $d$. This is because for
$\sigma _{\pm}\in \Gamma (L^{\pm 1/2})$, the following determines
$\nabla_{\frac{A_{0}}{2}}(\sigma_{-})$ :
 \begin{eqnarray*}
\nabla_{A_{0}} (\sigma_{+}+0)\otimes \sigma_{-}&=&
\partiali (\sigma_{+}+0 )\otimes \sigma_{-} + (\sigma_{+}+ 0 )\otimes
 \nabla_{\frac{A_{0}}{2}}(\sigma_{-})
\\ =\nabla_{A_{0}}(\sigma_{+}\otimes \sigma_{-})&=&
 d(\sigma_{+}\otimes \sigma_{-})
\end{eqnarray*}

The following is essentially the Leibnitz formula for Laplacian  applied to
Weitzenbock formula (7)

\vspace{.12in}

\noindent{\bf Proposition}:
 Let   $A, A_{0} \in \cal{A}(L^{-1}) $  and $i\;a=A-A_{0}$. Let
$\nabla _{a}=d +i\;a $ be the  covariant derivative of the trivial bundle
$\underline{\C}\longrightarrow X$, and $\alpha :X\to \C$. Let
$u_{0}$ be a section of $W^{+}=L^{-1}\oplus{\C}$ with  a constant
$\C $ component and $\Di_{A_{0}}(u_{0})=0$ then:
\begin{eqnarray}
\Di_{A}^{2}(\alpha u_{0})= (\nabla _{a}^{*}\nabla _{a} \alpha )u_{0} +
\frac{1}{2} \rho (F_{a} )\;\alpha \;u_{0} -2<\nabla _{a}\alpha \;,
\nabla_{A_{0}}(u_{0})>  \end{eqnarray}


Proof: By writing $\nabla_{A}=\nabla ^{A}$ for the sake of not cluttering
notations,
and abbreviating
$\;\nabla_{e_{j}}=\nabla_{j}\;$ and $\;\nabla^{a}_{j}(\alpha)=
\nabla_{j}(\alpha) + i\; a_{j} \alpha \;$, and leaving out summation signs for
repeated indices (Einstein convention) we calculate:

\begin{eqnarray}
\nabla^{A}(\alpha u_{0})&=&\nabla^{A}(\alpha)u_{0}+\alpha
\nabla^{A}(u_{0})\nonumber \\ &=& e^{j}\otimes\nabla_{j}(\alpha)u_{0}+
\alpha (\nabla^{A_{0}}(u_{0}) +i\; e^{j}\otimes a_{j}\;u_{0}) \nonumber
\nonumber \\ &=&
e^{j}\otimes(\nabla_{j}(\alpha) +i\; a_{j} \alpha\;) u_{0} + \alpha
\nabla^{A_{0}}(u_{0}) \nonumber \\
\Di_{A}(\alpha u_{0} )&=&\rho ( e^{j}) \;\nabla ^{a}_{j}(\alpha)\; u_{0} +
\alpha \;\Di_{A_{0}}(u_{0}) = \rho ( e^{j}) \;\nabla ^{a}_{j}(\alpha)\; u_{0}
\end{eqnarray}

\vspace{.1in}

\noindent By abbreviating $\;\mu=\nabla ^{a}_{j}(\alpha)\; $ we calculate:

\begin{eqnarray}
\nabla^{A} ( \rho ( e^{j}) \;\mu\; u_{0} )&=& e^{k}\otimes \rho ( e^{j})
\;\nabla_{k}(\mu )  u_{0} + e^{k}\otimes\; \rho(e^{j})\;\mu
\; (\nabla^{A_{0}}_{k}(u_{0}) +i\; a_{k}\;u_{0}) \nonumber \\ &=& e^{k}\otimes
\rho(e^{j}) \;\nabla^{a}_{k}(\mu)\; u_{0}\; + e^{k}\otimes
\rho(e^{j})\;\mu\;\nabla^{A_{0}}_{k}(u_{0})  \nonumber\\
\Di_{A}(  \rho ( e^{j}) \;\mu\; u_{0})&=& \rho (e^{k}
)\rho(e^{j})\;\nabla^{a}_{k}(\mu)\; u_{0}\; +
\rho (e^{k} )\rho(e^{j})\;\mu\;\nabla^{A_{0}}_{k}(u_{0}) \nonumber \\ & = &
-\nabla^{a}_{j}(\mu)\; u_{0}\; + \frac{1}{2}
\sum_{k,j} \rho (e^{k} )\rho(e^{j}) (\nabla^{a}_{k}(\mu)-
\nabla^{a}_{j}(\mu) )u_{0} \nonumber\\ & & - \mu\;\nabla^{A_{0}}_{j}(u_{0}) -
\mu\;\rho(e^{j})
\sum_{k\neq j} \rho (e^{k} ) \nabla^{A_{0}}_{k}(u_{0})
\end{eqnarray}
Since $0=\Di_{A_{0}}(u_{0})=\sum \rho (e^{k} )
\nabla^{A_{0}}_{k}(u_{0}) $ the last term of (3) is
$-\mu\;\nabla^{A_{0}}_{j}(u_{0}) $.

\vspace{.1in}

\noindent By plugging $\mu=\nabla ^{a}_{j}(\alpha)\; $ in (10) and summing over
$j$, from (2) we see
\begin{eqnarray*}
\Di_{A}^{2}(\alpha u_{0})=-\nabla ^{a}_{j}\nabla ^{a}_{j}(\alpha) u_{0}
+\frac{1}{2}
\rho(\sum F^{a}_{k,j}\; e^{k}\wedge e^{j})\;\alpha \;u_{0}-2\sum
\nabla ^{a}_{j}(\alpha) \nabla^{A_{0}}_{j}(u_{0})\;\;\;\;\Box
\end{eqnarray*}

\noindent{\bf Remark}: Notice that since $u_{0}$ has a constant
$\C$ component and
$\nabla_{A_{0}}$ restricts to the usual $d$ the $\C$ component, the term
$<\nabla _{a}\alpha \;, \nabla_{A_{0}}(u_{0})>$ lies entirely in $L$
component of
$W^{+}$

\section{Seiberg-Witten invariants} Let $X$ be a  closed oriented Riemanian
manifold, and $L\longrightarrow X$  a characteristic line bundle. Seiberg
-Witten
equations are defined for
$(A,\psi)\in {\cal A}(L)\times \Gamma(W^{+})$,
\begin{eqnarray}
\Di_{A}(\psi)&=&0\\
\rho(F_{A}^{+})&=&\sigma(\psi)
\end{eqnarray}
  Gauge group $\;{\cal G}(L)=Map(X,S^{1})\;$ acts on
$\;\tilde{\cal B}(L)={\cal A}(L)\times \Gamma(W^{+})\;$ as follows: for
$\;s=e^{if}\in {\cal G}(L) $
$$s^{*}(A, \psi)= (s^{*} A, s^{-1}\psi)=(A+s^{-1}ds\;,\; s^{-1}
\psi)= (A+i\; df \; ,\;s^{-1}\psi) $$

\noindent By locally writing $W^{\pm}= V^{\pm}\otimes L^{1/2} $, and
$\psi= \varphi\otimes \lambda\in \Gamma ( V^{\pm} \otimes L^{1/2})$ and from:
\begin{eqnarray*}
\Di _{s^{*}A}( \varphi\otimes \lambda)=\; \partiali (\varphi)\otimes
\lambda + [\; \varphi\otimes D_{A}(\lambda) + i\;df\; ( \varphi
\otimes \lambda)\;]\;
\end{eqnarray*} we see that
$\Di _{s^{*}A}\; (s^{-1}\psi) =s^{-1} \Di_{A} (\psi) $, and from definitions
$$\rho(F_{s^{*}A}^{+})=s^{-1}\rho(F_{A}^{+})\;s =
\rho(F_{A}^{+})=\sigma(\psi) =\sigma(s^{-1}\psi) $$
\noindent Hence the solution set  $\;\tilde{\cal M}(L)\subset
\tilde{\cal B}(L)\;$ of Seiberg-Witten equations is preserved by the action
$\;(A,\psi)\longmapsto s^{*}( A,\psi )\;$ of ${\cal G }(L)$ on
$\tilde{\cal M}(L)$. Define

$$ {\cal M}(L)=\tilde{\cal M}(L)/{\cal G }(L)  \;\;\subset\;\; {\cal B}(L)=
\tilde{\cal B}(L)/{\cal G} (L)$$

 We call a solution $(A,\psi)$  of (11) and (12) an irreducible solution if
$\psi\neq 0 $. $\;{\cal G}(L)$ acts on the subset  $\tilde{\cal M}^{*}(L)$
of the
irreducible solutions freely, we denote

$$ {\cal M}^{*}(L)=\tilde{\cal M}^{*}(L)/{\cal G }(L) $$
\vspace{.005in}

Any solution $(A,\psi) $ of Seiberg -Witten equations satisfies the
$C^{0}$ bound
\begin{eqnarray} |\psi|^{2}\leq \mbox{max}(0,-2s)
\end{eqnarray} where $s$ is the scalar curvature function of $X$. This
follows by
plugging (12) in the Weitzenbock formula (7).
  \begin{eqnarray}
\Di_{A}^{2}(\psi)&=& \nabla _{A}^{*}\nabla _{A}\psi
 +\frac{s}{4}\psi + \frac{1}{4}\sigma (\psi)\psi
\end{eqnarray} Then at the points where  where $|\psi|^{2}$ is maximum, we
calculate
\begin{eqnarray*} 0\leq \frac{1}{2}\Delta |\psi |^{2}&=& \frac{1}{2}
d^{*}d<\;\psi,\psi\;> =
\frac{1}{2} d^{*} (\; <\nabla_{A} \psi,\psi>+<\psi,\nabla_{A} \psi>
\;) \\ &=& \frac{1}{2} d^{*} ( \;\bar{<\psi, \nabla _{A} \psi>}+
<\psi,\nabla_{A}
\psi>\; )= d^{*} < \psi\;, \nabla_{A} \psi>_{\;\R}\\ &=& <\psi, \nabla_{A}
^{*}\nabla_{A} \psi>-|\nabla_{A}
\psi|^{2} \;\leq \; <\psi, \nabla_{A} ^{*}\nabla _{A} \psi>\\ &\leq & - \;
\frac{s}{4} |\; \psi \; |^{2} - \frac{1}{8} |\;\psi \; |^{4}
\end{eqnarray*} The last step follows from (14), (11) and (5), and the last
inequality gives (13)

\vspace{.15in}

\noindent{\bf Proposition 3.1} $\;{\cal M}(L)$ is compact

\vspace{.12in}

Proof: Given a sequence  $[\;A_{n},\psi_{n}\;] \in {\cal M}(L) $ we claim
that there
is a convergent subsequence (which we will denote by the same index), i.e.
there is
a sequence of gauge transformations $g_{n}\in {\cal G}(L)$ such that
$g_{n}^{*}(A_{n},\psi_{n})$  converges in $C^{\infty}$. Let $A_{0}$ be a base
connection. By Hodge theory of the elliptic complex:
$$\Omega^{0}(X)\stackrel {d^{0}}{\longrightarrow} \Omega^{1}(X)
\stackrel{d^{+}}{\longrightarrow} \Omega^{2}_{+}(X) $$
$$A-A_{0}=h_{n}+a_{n}+b_{n}\in {\cal H}\oplus im (d^{+})^{*}\oplus im (d)
$$ where
$\cal{H}$ are the  harmonic $1$-forms. After applying gauge transformation
$g_{n}$
we can assume that $b_{n}=0$, i.e. if
$b_{n}=i\;d f _{n}$ we can let
$g_{n}=e^{if}$. Also $$h_{n}\in {\cal H}=H^{1}(X;\R)\;\;\;\mbox{and a
component of}\;\;{\cal G}(L)\;
\;\mbox{is}\;\; H^{1}(X;\Z)$$ Hence after a gauge transformation we can assume
$ h_{n}\in H^{1}(X;\R)/ H^{1}(X;\Z)$ so $h_{n}$ has convergent
subsequence. Consider the first order elliptic operator:

$$D =d^{*} \oplus d^{+} :
\Omega^{1}(X)_{L^{p}_{k}}\longrightarrow \Omega^{0}(X)_{L^{p}_{k-1}}
\oplus \Omega^{2}_{+}(X)_{L^{p}_{k-1}}$$

\noindent The kernel of $D$ consists of harmonic $1$-forms, hence by Poincare
inequality
 if $a$ is a $1$-form orthogonal to the harmonic forms, then for some
constant $C$
$$ ||a ||_{L^{p}_{k}}\leq C ||D(a)||_{L^{p}_{k-1}} $$

\noindent Now  $a_{n}=(d^{+})^{*} \alpha_{n}$ implies
$d^{*}(a_{n})=0 $. Since $\alpha_{n}$ is orthogonal to harmonic forms, and
by calling
$A_{n}=A_{0}+a_{n}$ we see :

$$ || a_{n}||_{L^{p}_{1}}\leq C\; ||D(a_{n}) ||_{L^{p}}
\leq C || d^{+} a_{n} ||_{L^{p}}= C \;|| F_{A_{n}}^{+}-F_{A_{0}}^{+}
||_{L^{p}}$$
Here  we use C for a generic constant. By (12), (4) and (13) there is a $C$
depending only on the scalar curvature $s$ with
\begin{eqnarray}|| a_{n}||_{L^{p}_{1}}\leq C \end{eqnarray} By iterating this
process we get $ || a_{n}||_{L^{p}_{k}}\leq C $ for all
$k$ , hence
$|| a_{n}||_{\infty}\leq C$. From the elliptic estimate and
$\Di_{A_{n}}(\psi_{n})=0$ :

\begin{eqnarray} ||\psi_{n}||_{L^{p}_{1}}&\leq & C(\; ||\Di
_{A_{0}}\psi_{n}\;||_{L^{p}} + || \psi _{n} ||_{L^{p}})= C(\;
||a_{n}\psi_{n}\;||_{L^{p}}+|| \psi _{n} ||_{L^{p}}) \nonumber \\
||\psi_{n}||_{L^{p}_{1}}
&\leq &    C(\; ||a_{n}||_{\infty} ||\psi _{n} ||_{L^{p}} + ||\psi _{n}
||_{L^{p}} )
\leq C
\end{eqnarray}

By repeating this (boothstrapping) process  we get
$||\psi_{n}||_{L^{p}_{k}}\leq  C $, for all $k$, where C depends only on
the scalar
curvature $s$ and $A_{0}$. By Rallich theorem we get convergent subsequence of
$\;(a_{n},\psi_{n})\;$ in
${L^{p}_{k-1}} $ norm for all $k$. So we get this convergence to be
$C^{\infty}$ convergence.
\hspace{2in}$\Box$

\vspace{.15in} It is not clear that the solution set of  Seiberg-Witten
equations is
a smooth manifold. However we can  perturb the Seiberg-Witten equations
(11), (12)
by any self dual
$2$-form $\delta \in \Omega^{+}(X)$, in a gauge invariant way, to obtain a
new set
of equations whose solutions set is a smooth manifold:
\begin{eqnarray}
\Di_{A}(\psi)&=&0\\
\rho(F_{A}^{+} + i\;\delta )&=&\sigma(\psi)
\end{eqnarray}
\vspace{.15in} Denote this solution space by $\tilde{\cal M}_{\delta}(L)$, and
parametrized solution space by
$$\tilde{\bf{\cal M}}=\bigcup_{\delta \in \Omega^{+}}
\tilde{\cal M}_{\delta} (L)\times \{\;\delta\;\}\subset {\cal A}(L)\times
\Gamma(W^{+})
\times \Omega ^{+}(X)$$
 $${\cal M}_{\delta}(L) = \tilde{\cal M}_{\delta}(L)\;/{\cal G}(L)
\;\;\subset\;\; {\bf{\cal M}} =\tilde{\bf{\cal M}}\;/{\cal G}(L) $$

\vspace{.08in}

\noindent Let $\tilde{\cal M}_{\delta}(L)^{*}\subset\tilde{\bf{\cal
M}}^{*}$ be the
corresponding irreducible solutions, and also let
$ {\cal M}_{\delta}(L)^{*} \subset {\bf{\cal M}}^{*}$ be their quotients by
Gauge
group. The following theorem says that for a generic choice of $\delta $
the set
${\cal M}_{\delta}(L)^{*}$ is a closed smooth manifold.

\vspace{.15in}
\noindent {\bf Proposition 3.2} $\;{\cal M}^{*}$ is a smooth manifold.
Projection
$\pi:{\cal M}^{*}\longrightarrow \Omega^{+}(X) $ is a proper surjection of
Fredholm
index:
$$d(L)=\frac{1}{4}[\;c_{1}(L)^{2}-(2\chi +3\sigma)\;]$$ where $\chi$ and
$\sigma$
are Euler characteristic and the signature of $X$.

\vspace{.12in}

Proof: The linearization of the map
$\;(A,\psi,\delta)\longmapsto (\rho(F_{A}^{+} + i\;\delta )-
\sigma(\psi) ,\Di_{A}(\psi)\;) $ at $(A_{0},\psi_{0},\delta_{0})$ is given by:
$$ P: \Omega^{1}(X)\oplus \Gamma(W^{+})\oplus
\Omega^{+}(X)\longrightarrow su(W^{+})\oplus \Gamma(W^{-})$$
$$ P(a,\psi, \epsilon)=(\rho (d^{+}a + i\; \epsilon ) -2
\;\sigma(\psi,\psi_{0})\;,\;
\Di_{A_{0}}\psi +\rho(a)\psi_{0})  $$ To see that this is onto we pick $
(\kappa,\theta)\in su(W^{+})\oplus
\Gamma(W^{-})$, by varying $\epsilon$ we can see that $(\kappa,0) $ is in
the image
of $P$. To see
$(0,\theta)$ is in the image of $P$, we prove that if it is in the orthogonal
complement to
$\;image(P)\;$ then it is $(0,0)$; i.e. assume
$$<\Di_{A_{0}}\psi, \theta> +<\rho(a)\psi_{0}, \theta >=0$$ for all
$a$ and $\psi$, then by choosing $\psi=0$ we see $<\rho(a)\psi_{0},
\theta >=0$ for all $a$ which implies $\theta=0$

\vspace{.05in}

By implicit function theorem $\;\tilde {\cal M}\;$ is a smooth manifold, and by
Sard's theorem $\;\tilde{\cal M}_{\delta}(L)\;$ are smooth manifolds, for
generic
choice of $\delta $'s. Hence their free quotients
$\;{\bf{\cal M}}^{*\;}$  and $\;{\cal M}_{\delta}(L)^{*}\;$ are smooth
manifolds.

After taking ``gauge fixig" account, the dimension of $ {\cal M}_{\delta}(L)$
is
given by the index of
$\;P+ d^{*}$ (c.f. [{\bf DK}]).  $\;P+ d^{*}$ is the compact perturbation of

$$ S: \Omega^{1}(X)\oplus
\Gamma(W^{+})\longrightarrow [\;\Omega^{0}(X)\oplus
\Omega_{+}^{2}(X)\;]\oplus \Gamma(W^{-}) $$

$$S=\left( \begin{array}{cc} d^{*}\oplus d^{+} & 0 \\ 0 & \Di_{A_{0}}
\end{array} \right)$$

By Atiyah-Singer index theorem
\begin{eqnarray}
\mbox{ dim }{\cal M}_{\delta} (L)=\mbox{ind} (S)&=&
 \mbox{index} (d^{*}\oplus d^{+}) + \mbox{index}
_{\R}\;\Di_{A_{0}}\nonumber\\ &=&
-\frac{1}{2} (\chi + \sigma )+\frac{1}{4}( c_{1}(L)^{2}-\sigma ) \nonumber
\\ &=&
\frac{1}{4}\;[\;c_{1}(L)^{2}-(2\chi + 3 \sigma )\;]\nonumber \\ &=&
\frac{c_{1}(L)^{2}-\sigma
}{4}\;-\;(1+b^{+})
\end{eqnarray} where $b^{+}$ is the dimension of positive define part
$\;H^{2}_{+}\;$ of
$\;H^{2}(X;\Z)$. Notice that when $b^{+}$ is odd this expression is even,
since $L$ being a characteristic line bundle we have
$c_{1}(L)^{2}=\sigma\; \mbox{mod}\; 8 $
\hspace{3in} $\Box$

\vspace{.15in}

Now assume that $H^{1}(X)=0$, then ${\cal G}(L)=K(\Z,1)$. Than being a free
quotient of a contractible space by ${\cal G}(L)$ we have
$${\cal B}^{*}(L)=K(\Z,2)=\C\P^{\infty}$$ The orientation of
$H^{2}_{+}$ gives an orientation to
${\cal M}_{\delta}(L)$. Now  By (19) if $ b^{+} $ is odd
 ${\cal M}_{\delta}(L)\subset{\cal B}^{*}(L) $ is an even dimensional $2d$
 smooth closed oriented submanifold, then we can define Seiberg-Witten
invariants as:
$$SW_{L}(X)=<{\cal M}_{\delta}(L)\;,\; [\;\C\P^{d}\;]> $$

As in the case of Donaldson invariants ([{\bf DK}]), even though
${\cal M}_{\delta}(L)$ depends on metric (and on the perturbation
$\delta$) the invariant $SW_{L}(X)$ is independent of these choices, provided
$\;b^{+}\ge 2\;$, i.e. there is a generic metric theorem.

Also by (13) if
$X$ has
nonnegative scalar curvature then all the solutions are reducible, i.e.
$\psi=0$.
This implies that $A$ is anti-self-dual, i.e. $F_{A}^{+}=0$; but just as in
[{\bf
DK}] , If
$b^{+}\geq 2$ and $L$ nontrivial, for a generic metric $L$ can not admit such
connections. Hence $\tilde{\cal M}=\emptyset$ which implies $SW_{L}(X)=0$.

Similar to Donaldson invariants there is a ``connected sum theorem" for
Seiberg-Witten invariants: If $X_{i}\;i=1,2\;$ are oriented compact smooth
manifolds
with common boundary, which is a
$3$-manifold with a positive scalar curvature; then gluing these manifolds
together
along their  boundaries produces a  manifold
$X=X_{1}\smile X_{2}$ with vanishing Seiberg-Witten invariants  (cf [{\bf
F}],[{\bf
FS}]). There is also conjecture that only
$0$-dimensional moduli spaces ${\cal M}_{\delta}(L)$ give nonzero invariants
$SW_{L}(X)$.

\section{Almost Complex and Symplectic Structures}

Now assume that $X$ has an almost complex structure. This means that there is a
principal
$GL(2,\C)$-bundle $Q\longrightarrow X$ such that
$$T(X)\cong Q\times_{GL(2,\C)} \C^{2}$$ By choosing Hermitian metric on
$T(X)$ we can assume $Q\longrightarrow X$ is a
$U(2)$ bundle, and the tangent frame bundle $P_{SO(4)}(TX)$ comes from $Q$
by the
reduction map
$$U(2)=(\;S^{1}\times SU (2) \;)/\Z_{2}\hookrightarrow (\;SU (2)\times SU
(2)\;)/\Z_{2} =SO(4)$$ Equivalently there is an endemorphism  $I\in
\Gamma(End(TX))$ with $I^{2}=-Id$
$$\begin{array}{ccc} T(X)& \stackrel{I}{\longrightarrow} & T(X)\\

\;\;\;\;\;\;\;\;\;\;\;\;\;\;\;\searrow & &
\swarrow\;\;\;\;\;\;\;\;\;\;\;\;\;\;\;\\ & X &
\end{array}$$ The $\;\pm i \;$  eigenspaces of $I$ splits the complexified
tangent
space
$T(X)_{\C}$
 $$T(X)_{\C}\cong T^{1,0}(X)\oplus T^{0,1}(X)=\Lambda^{1,0}(X)\oplus
\Lambda^{0,1}(X)$$ This gives us a complex line bundle which is called the
canonical
line bundle:
$$K=K_{X}=\Lambda^{2,0}(X)=\Lambda^{2}(T^{1,0})\longrightarrow X$$

Both $K^{\pm}$ are characteristic; corresponding to line bundle $K
\longrightarrow X$ there is a canonical $Spin_{c}(4)$ structure on
$X$, given by the lifting of $f[\lambda ,  A ]=([\lambda ,A],
\lambda^{2})$
$$\begin{array}{ccc} && Spin_{c}(4)\\ & & \\ & F \nearrow & l\downarrow\\ & &
\\
U(2) &\stackrel{f}{\longrightarrow} & SO(4)\times S^{1}
\end{array}$$
  $\;F[\lambda ,A]=[\lambda  , A,\lambda]$. Transition function
$\lambda^{2}$ gives $K$, and the corresponding
${\C}^{2}$-bundles are given by:
\begin{eqnarray*} W^{+}&= & \Lambda ^{0,2}(X)\oplus\Lambda^{0,0}(X) =
K^{-1}\oplus{\C}\\ W^{-}&=&\Lambda^{0,1}(X)
\end{eqnarray*} We can check  this from the transition functions, e.g. for
$W^{+}$,
$x=z+jw\in{\H}$

$$x\longmapsto \lambda x \lambda^{-1} =\lambda (z+jw) \bar{\lambda } = z + jw
\bar{\lambda} \bar{ \lambda} = z +jw \lambda ^{-2}$$

Since we can identify $\bar{\Lambda}^{0,1}(X)\cong
\Lambda^{1,0}(X)$, and
$\Lambda^{0,2}(X)\otimes \Lambda^{1,0}(X) \cong \Lambda^{0,1}(X)$ we
readily see the
decomposition
$T(X)_{\C}\cong W^{+}\otimes \bar{W}^{-}$. As real bundles we have

$$ \Lambda^{+}(X)\cong K\oplus {\R}$$

\noindent We can verify this by taking
 $\;\{e^{1},e^{2}=I(e^{1}),e^{3},e^{4}=I(e^{3} )\}\;$ to be a local orthonormal
basis for $T^{*}(X)$, then
$$\Lambda^{1,0}(X)=\;<e^{1}-ie^{2},e^{3}-ie^{4}>\;\;\;,and\;\;\;\;
\Lambda^{0,1}(X)=\;<e^{1}+ie^{2},e^{3}+ie^{4}>$$
$$K =\;< f =(e^{1}-ie^{2})\wedge (e^{3}-ie^{4})>$$
$$\Lambda^{+}(X)=\;<\omega=\frac{1}{2}( e^{1}\wedge e^{2} + e^{3}\wedge
e^{4}),\;f_{2}=
 \frac{1}{2}( e^{1}\wedge e^{3} + e^{4}\wedge e^{2}),\; f_{3}=
 \frac{1}{2}( e^{1}\wedge e^{4} + e^{2}\wedge e^{3})>  $$
$\omega $ is the global form $ \omega(X,Y)=g(X,IY) $ where $ g $ is the
hermitian
metric (which makes the basis
$\{e^{1},e^{2},e^{3},e^{4}\} $ orthogonal). Also since
 $f= 2(f_{2} -i f_{3}) $, we see as ${\R}^{3}$-bundles
$\Lambda^{+}(X)\cong K\oplus {\R}(\omega)$. We can check:
$$W^{+}\otimes \bar{W}^{+}\cong{\C}\oplus{\C}\oplus K\oplus
\bar{K}= (K\oplus{\R})_{\C}\oplus {\C}$$
 As before by writing the sections of $W^{+}$ by $z+jw \in \Gamma ({\C}\oplus
K^{-1})$
 we see that $\omega, f_{2},f_{3}$ act as Pauli matrices; in particular
\begin{eqnarray*}
\omega & \longmapsto & {\left(\begin{array}{cc} i &0 \\ 0 &- i
\end{array}\right)}\\ f &\longmapsto &{2\left(\begin{array}{cc} 0 & -1 \\ 1 & 0
\end{array}\right)- 2i\left(\begin{array}{cc} 0 & -i
\\- i & 0  \end{array}\right)=
\left(\begin{array}{cc} 0 & -4 \\ 0 & 0  \end{array}\right)}\\
\bar{f} &\longmapsto &{2\left(\begin{array}{cc} 0 & -1 \\ 1 & 0
\end{array}\right)+ 2i\left(\begin{array}{cc} 0 &- i \\ -i & 0
\end{array}\right)=
\left(\begin{array}{cc} 0 & 0 \\ 4 & 0  \end{array}\right)}\\
\end{eqnarray*}

So in particular, if we write $\psi\in
\Gamma (W^{+})=\Gamma({\C}\oplus K^{-1})$ by $\psi=\alpha u_{0} + {\beta} \;$,
where $\beta$ is a section of $K^{-1}$, and
$\alpha: X\to {\C} $ and $u_{0}$ is a fixed section of the trivial bundle
$\underline{\C}\to X$ with $||u_{0}||=1$, then
$$\rho(\omega  )\; u_{0}= iu_{0} \;\;\;\; \rho(\omega
)\;\beta=-i\beta\;\;\;\;\;\;\;\;\;\;\;\;\;\;\;\;\;$$
$$\rho(\beta)\; u_{0}=4\beta\;\;\;\; \rho(\beta
)\;\beta=0\;\;\;\;\;\;\;\;\;\;\;\;\;\;\;\;(*)$$
$$\;\;\;\;\;\;\rho(\bar{\beta})\; u_{0}=0
\;\;\;\;\;\;\;\;\;\rho(\bar{\beta} )
\;\beta =-4\;|\beta |^{2}u_{0}
\;\;\;\;\;\;\;\;\;\;\;\;\;\;\;\;\;\;$$

\noindent We see these  by locally  writing $\psi$ in terms of basis
$\;\psi =  \alpha u_{0} + \lambda \bar{f} $,  where $\beta= \lambda
\bar{f}$ with
$\;||\bar{f}||=1 $. Writing Formula (3) in terms of the basis $\{
\omega, f, \bar{f}\}$ we get:
\begin{eqnarray}i\; \sigma ( \alpha, \lambda)& =&
\rho\;[\;\frac{|\alpha |^{2}-|\lambda |^{2}}{2}\;\omega -
\frac{i}{4}\;\alpha  \bar {\lambda} f +\frac{i}{4}
\;\bar{\alpha}{\lambda} \bar{f}\;] \nonumber\\
\sigma (\psi)&=& \rho\;[\;\frac{|\beta |^{2}-|\alpha |^{2}}{2}\;i\;
\omega - \frac{1}{4}\;\alpha \; \bar {\beta} +\frac{1}{4}
\;\bar{\alpha}\; \beta \;]
\end{eqnarray}

If we consider the decomposition $\;F_{A}^{+}=F^{2,0}_{A} +F^{0,2}_{A}
+F^{1,1}_{A}\; $
 the equation $\rho(F_{A})=\sigma(\psi)$ gives  Witten's formulas:

\begin{eqnarray} F^{2,0}_{A}&=& -\frac{1}{4}\alpha \;\bar{\beta}\\
F^{0,2}_{A}&=&
\;\frac{1}{4}\bar{\alpha} \;\beta\\ F^{1,1}_{A}&=&
\frac{|\beta |^{2}-|\alpha |^{2}}{2}\;i\;\omega
\end{eqnarray}

\vspace{.18in}

In case $\;X\;$ is a  Kahler surface the Dirac operator is given by
(c.f.[{\bf LM}])
$$\Di_{A}=\bar{\partiali}^{*}_{A} + \bar{\partiali}_{A} :\Gamma (W^{+})\to
\Gamma (W^{-})$$ Hence from the Dirac part of the Seiberg-Witten equation
(17) we
have
\begin{eqnarray}
\bar{\partiali}_{A}^{*}(\beta) + \bar{\partiali}_{A}(\alpha u_{0})&=&0
\nonumber \\
\bar{\partiali}_{A} \bar{\partiali}^{*}_{A}(\beta) +
\bar{\partiali}_{A}\bar{\partiali}_{A}(\alpha u_{0})&=&0
\end{eqnarray}

\noindent The second term is $\;
\bar{\partiali}_{A}\bar{\partiali}_{A}(\alpha u_{0})= F^{0,2}_{A}
\alpha u_{0}= \frac{1}{4} |\alpha|^{2}\beta $.  By taking inner product
both sides
of (24) by $\;\beta \;$ and integrating over $X$ we get the
$L^{2}$ norms satisfy
\begin{eqnarray}||\alpha||^{2} ||\beta||^{2}=0\;\;& \Longrightarrow
\;\;\alpha= 0\; \mbox{ or }\;\beta=0
\end{eqnarray} This argument  eventually calculates
$\;SW_{K}(X)=1\;$ ([{\bf W}]). We will not repeat this argument here, instead
we
will review a stronger result of C.Taubes for symplectic manifolds below, which
implies this result.

\vspace{.18in}

We call an almost complex manifold with Hermitian metric $\{X, I,g\}$
syplectic if
$d\omega=0$. Clearly a nondegenerate closed form $\omega$ and a hermitian
metric
determines the almost complex structure $I$. Given
$\omega$ then $I$ is called an almost complex structure taming the
symplectic form
$\omega$

\vspace{.1in}

By Section 2.1 there is a unique connection $A_{0}$ in $\; K\longrightarrow
X \;$
such that the induced Dirac operator
$D_{A_{0}}$ on $W^{+}$ restricted to the trivial summand $\underline {\C}\to X$
is the exterior derivative $d$. Let $u_{0}$ be the section of $W^{+}$ with
constant
${\C}$ component and
$||u_{0}||=1$. Taubs's first fundamental observation is
$$\Di_{A} (u_{0})=0 \;\;\;\; \mbox{ if and only if}\;\;\;\;\; d\omega=0 $$
This can
be seen by applying the Dirac operator to both sides of $ i u_{0} =
\rho(\omega).
u_{0}$, and observing that by the choice of $u_{0}$ the term
$\nabla_{A_{0}}(u_{0})$
lies entirely in
$K^{-1}$ component:
\begin{eqnarray*} i \Di _{A_{0}}(u_{0})&=&\sum
\rho(e^{i})\nabla_{i}\;(\rho (\omega )\;u_{0}) \\ &=&\sum \rho (e^{i}) \;[\
\nabla_{i}\;(\rho(\omega) )\;u_{0} +
\rho(\omega)\;\nabla_{i}\;(u_{0})\;]\\ &=& \sum \rho (e^{i})
\nabla_{i}\;(\rho(\omega) ) \;u_{0} -i\;\sum\rho (e^{i}\;)\nabla
_{i}\;(u_{0})\\
2i\;\Di _{A_{0}}(u_{0}) &=& \sum \rho (e^{i})
\nabla_{i}\;(\rho(\omega) )
\;u_{0}=
\rho ( (d +d^{*}) \omega)\;u_{0}=\rho ((d-*d)\omega)\;u_{0}
\end{eqnarray*} Last equality holds since $\omega \in
\Lambda^{+}(X)_{\C}\oplus {\C}$, and by naturality, the Dirac operator on
$ \Lambda^{*}(X)_{\C} $ is
$\;d+d^{*}$, and  since $d=-*d*$ on $2$ forms and $\omega$ is self dual
$$2i\;\Di _{A_{0}}(u_{0})=-\rho (*d \omega ) u_{0}$$

\vspace{.15in}

\noindent{\bf Theorem (Taubes) }: Let $(X,\omega)$ be a closed symplectic
manifold
such that $b_{2}(X)^{+}\geq2$, then
$SW_{K}(X)=\pm1 $.

\vspace{.1in}

Proof: $\mbox{Write}\;\;\;\psi=\alpha u_{0} + {\beta}  \in \Gamma( W^{+})=
\Gamma ({\C}\oplus K^{-1}) $ where $\alpha: X\to {\C} $, and
$u_{0}$ is the section as above. Consider the perturbed Seiberg-Witten
equations :
For
$(A,\psi)\in {\cal A}(L)\times \Gamma (W^{+})$ :
\begin{eqnarray} \Di _{A}(\psi)&=&0\\
\rho (F_{A}^{+})&=&\rho (F_{A_{0}}^{+})+r\;[\;
\sigma(\psi)+i\;\rho(\omega)\;]
\end{eqnarray} By (20) the second equation is equivalent to:
\begin{eqnarray}
 F_{A}^{+}-F_{A_{0}}^{+} &=&r\;\ [\; ( \;\frac{|\beta
|^{2}-|\alpha|^{2}}{2} +1 )\;
i \omega - \frac{1}{4}\alpha
\bar{\beta} +\frac{1}{4} \bar{\alpha} \beta \; ]
\end{eqnarray}

We will show that up to gauge equivalence there is a unique solution to these
equations.
 Write $A=A_{0}+a$, after a gauge transformation we can assume that
$a$ is coclosed, i.e.
$d^{*}(a)=0$. Clearly $ (A,\psi)=(A_{0},u_{0})$, and $r=0$ satisfy these
equations.
It suffices to show that for $r\longmapsto \infty $ these  equations admit only
$(A_{0},u_{0})$ as a solution. From Weitzenbock formulas (7), (8) and
abbreviating
$\nabla_{A_{0}}(u_{0})=b$ we get

\begin{eqnarray}
\Di_{A}^{2}(\psi)= \Di_{A}^{2}(\beta)+(\nabla _{a}^{*}\nabla_{a}\alpha)u_{0}
 -2<\nabla_{a} \alpha,b> + \frac{1}{2}\alpha\;
\rho(F_{A}^{+}-F_{A_{0}}^{+})\;u_{0}\\
\Di_{A}^{2}(\beta) =(\nabla _{A}^{*}\nabla_{A}\;\beta) +\frac{s}{4}\;\beta
+\frac{1}{4} \rho ( F_{A_{0}}^{+})\;\beta +\frac{1}{4} \rho
(F_{A}^{+}-F_{A_{0}}^{+})
\;\beta
\hspace{.3in}
\end{eqnarray}

\noindent From (28) and (*) we see that
\begin{eqnarray}\frac{1}{2}\alpha
\;\rho(F_{A}^{+}-F_{A_{0}}^{+})u_{0}&=&
\frac {r}{4}\alpha \;(\;|\alpha|^{2}-|\beta|^{2}-2)\;u_{0} +
\frac{r}{2} \;|\alpha |^{2}\beta \\
\frac{1}{4}\;\rho(F_{A}^{+}-F_{A_{0}}^{+})\;\beta &= & -
\;\frac{r}{8}\;(\;|\alpha|^{2}-|\beta|^{2}-2)
\;\beta +\frac{r}{4}\alpha |\beta|^{2} u_{0}
 \end{eqnarray}

\noindent By substituting (31) in (29) we get

\begin{eqnarray}\Di_{A}^{2}(\psi -\beta)&=& [\;\nabla _{a}^{*}\nabla_{a}\alpha
+
\frac{r}{4}\alpha
\;(|\alpha|^{2}-|\beta|^{2}-2)\;]\;u_{0}\nonumber\\ &&  -2 <\nabla_{a}
\alpha,b> +\frac{r}{2}\;
|\alpha\;|^{2} \beta
\end{eqnarray}

\noindent By substituting  (32) in (30), then substituting (30) in (33)  we
obtain:
\begin{eqnarray}0=\Di_{A}^{2}(\psi)&=& [\;\nabla _{a}^{*}\nabla_{a}\alpha +
\frac{r}{4}\alpha
\;(|\alpha|^{2}-2)\;]
\;u_{0}-2 <\nabla_{a} \alpha,b>  \nonumber \\ & &
+[\;\nabla _{A}^{*}\nabla_{A}+\frac{s}{4} +
\frac{1}{4}\;\rho(F_{A_{0}}^{+})+\frac{r}{8} (\;3|\alpha|^{2}+|\beta
|^{2}+2)\;]\;\beta
\end{eqnarray}

\noindent By recalling  that $\beta$ and $u_{0}$ are orthogonal sections of
$W^{+}$, we take inner product of both sides of (8) with $\beta$ and
integrate over
$X$ and obtain:
\begin{eqnarray*}
\int_{X}(\;|\nabla_{A}\beta\;|^{2} +\frac{r}{8}\;|\beta|^{4} +
\frac{r}{4}|\beta|^{2}+\frac{3r}{8}|\alpha|^{2}|\beta|^{2}\;)&=& \\
2\int_{X}(<<\nabla_{a} \alpha \;,\;b>\;,\;\beta> -\;\frac{s}{4}\;|\beta|^{2}
-\;\frac{1}{4}<\rho(F_{A_{0}}^{+})\beta
\;,\; \beta> & &
\end{eqnarray*}

$$\mbox{Hence}\;\; \int_{X} \;|\nabla_{A}\beta\;|^{2}
+\frac{r}{8}\;|\beta|^{4} +
\frac{r}{4}|\beta|^{2}+\frac{3r}{8}|\alpha|^{2}|\beta|^{2}\;
\leq  \int_{X}c_{1}|\beta|^{2}+ c_{2}|\beta | |\nabla_{a} \alpha|  $$ where
$c_1$
and $c_2$ are positive constants depending on the Riemanian metric and the base
connection $A_{0}$. Choose $r\gg 1$, by calling  $c_{2}=2c_{3}$ we get :

$$ \int_{X}(\;|\nabla_{A}\beta\;|^{2} +\frac{r}{8}\;|\beta|^{4} +
\frac{r}{8}|\beta|^{2}+\frac{3r}{8}|\alpha|^{2}|\beta|^{2}\;)
\leq  \int_{X}(c_{1}-\frac{r}{8})\;|\beta\;|^{2}+ 2c_{3}| \beta\; | |\nabla_{a}
\alpha|  $$
$$=- \left[ \;(r/8-c_{1})^{1/2}\;|\beta\;|-c_{3}(r/8-c_{1})^{-1/2}\;
|\nabla_{a}
\alpha|\;\right]^{2}+ \frac{c_{3}^{2}}{(r/8-c_{1})} |\nabla_{a} \alpha|^{2}
\leq \frac{C}{r}|\nabla_{a} \alpha|^{2}$$ For some $C$ depending on the
metric and
$A_{0}$. In particular we have
\begin{eqnarray*}
\int_{X}r \;|\beta|^{2} -\frac{8C}{r}\;|\nabla_{a} \alpha|^{2}&\leq &0 \\
8c_{2}\;|\beta \;|\;|\nabla_{a} \alpha|- \frac{8C}{r}|\nabla_{a}
\alpha|^{2}&\leq&
\int_{X}({r}-8c_{1})\;|\beta\;|^{2}
\end{eqnarray*}
\begin{eqnarray}\mbox{Hence}\;\;\;\;\;\;\;\; c_{2}\;|\beta
\;|\;|\nabla_{a} \alpha|- \frac{2C}{r}|\nabla_{a}
\alpha|^{2}&\leq &0
\end{eqnarray}

\vspace{.1in}

Now by self adjointness of the Dirac operator, and by $\alpha
u_{0}=\psi-\beta $ we
get:
\begin{eqnarray} <\Di_{A}^{2}(\psi)\;,\;\alpha u_{0}>&=&<\Di_{A}^{2}(\psi
-\beta)\;,\;\alpha u_{0}>+ <\Di_{A}^{2}\;(\beta)\;,\;\alpha u_{0}>\nonumber\\
&=&<\Di_{A}^{2}(\psi -\beta)\;,\;\alpha u_{0}>+
<\beta\;,\;\Di_{A}^{2}(\psi-\beta)>
\end{eqnarray}

\noindent We can calculate (36) by using (33) and obtain the inequalities:

\begin{eqnarray*} 0=<\Di_{A}^{2}(\psi)\;,\;\alpha
u_{0}>&=&|\nabla_{a}\alpha|^{2} +
\frac{r}{4}|\alpha|^{4}-
\frac{r}{4}|\alpha|^{2}|\beta|^{2}-
\frac{r}{2}|\alpha|^{2}\\ && +
\;\frac{r}{2}\;|\alpha|^{2}|\beta|^{2}-2 <<\nabla_{a} \alpha,b>,
\beta>
\end{eqnarray*}
\begin{eqnarray*}
\int_{X}|\nabla_{a}\alpha|^{2} +\frac{r}{4}|\alpha|^{4}-\frac{r}{2}|\alpha|^{2}
&\leq &\int_{X}2 <<\nabla_{a} \alpha,b>, \beta> -\frac{r}{4}
|\alpha|^{2}|\beta|^{2}\\ &\leq& \int_{X}2 <<\nabla_{a} \alpha,b>, \beta>\leq
\int_{X}c_{2} |\nabla_{a} \alpha|\; |\beta\; |
\end{eqnarray*} By choosing $c_{4}=1-{2C}/{r}$ and by (35), we see
\begin{eqnarray}
\int_{X}c_{4}\;|\nabla_{a}\alpha|^{2}
+\frac{r}{4}|\alpha|^{4}-\frac{r}{2}|\alpha|^{2} &\leq &0
\end{eqnarray} Since for a connection $A$ in $K\longrightarrow  X$ the class
$({i}/{2\pi}) F_{A}$ represents the Chern class
$c_{1}(K)$, and since $\omega$ is a self dual two form we can write:
$$\int_{X}\omega \wedge F_{A}=-2\pi\;i\; \omega c_{1}(K) \;\;\;\;\;\;
\int_{X}\omega \wedge F_{A}=\int_{X}\omega \wedge F_{A}^{+}$$
$$\int_{X} \omega \wedge  (F_{A}^{+}-F_{A_{0}}^{+})=0$$ By (28) this implies:
\begin{eqnarray}
\frac{r}{2} \int_{X} (2-|\alpha|^{2}+|\beta|^{2})=0
\end{eqnarray}

By adding (38) to (37)  we get

\begin{eqnarray}
 \int_{X}c_{4}\;|\nabla_{a}\alpha|^{2} + \frac{r}{2}|\beta|^{2} +r (1 -
\frac{1}{2}|\alpha|^{2})^{2}\leq 0
\end{eqnarray}
  Assume $r\gg 1$, then $c_{4}\geq 0$ and hence $\nabla_{a}\alpha=0$ and
$\beta=0$ and $|\alpha|=\sqrt{2}$, hence:
$$\beta =0\;\;\mbox{and}\;\;\alpha =\sqrt{2}e^{i\;\theta}\;\;\mbox{ and}\;\;
\nabla_{a}(\alpha)=d(e^{i\;\theta})+i\;a\;e^{i\;\theta}=0$$ Hence
$a= d(-\theta)$, recall that we also have $d^{*}(a)=0$ which gives
$$0 = <d^{*}d(\theta),\theta> = <d(\theta),d(\theta>)>=||d(\theta)||^{2}
$$ Hence
$\;a=0\;$  and
$\;\alpha \;$=constant. So up to a gauge equivalence
$(A,\psi)=(A_{0},u_{0})$

\hspace{4.5in} $\Box$

\section{Applications}

Let $X$ be a  simply connected closed smooth $4$-manifold. By
J.H.C.Whitehead the
intersection form
$$ q_{X}:H_{2}(X;\Z)\otimes H_{2}(X;\Z) \longrightarrow \Z$$
 determines the homotopy type of $X$. By C.T.C Wall in fact $q_{X}$
determines the
$h$-cobordism class of $X$. Donaldson (c.f. [{\bf DK}]) showed that if
$q_{X}$ is
definite then it is dioganalizable, i.e.
$$q_{X}=p<1>\oplus q<-1>$$ We call $q_{X}$ is even if q(a,a) is even for
all $a$,
otherwise we call
$q_{X}$ odd. Since  integral liftings $c$ of the second Steifel Whitney calass
$w_{2}$ of $X$ are characterized by $\;c.a=a.a\;$ for all $a\in H_{2}(X;\Z)$,
the condition of $q_{X}$ being even is equivalent to $X$ being spin. From
classification of unimodular even integral quadratic forms and the Rohlin
theorem it
follows that the intersection form of a closed smooth spin manifold is in
the form:
\begin{eqnarray}q_{X}=2k E_{8}\oplus lH  \end{eqnarray} where
$E_{8}$ is the $8\times 8$ intersection matrix given by the Dynikin diagram

\begin{figure}[htb]
\vspace{.8in}
\special{picture 1 scaled 700}
\caption{}
\end{figure}
\noindent and $H$ is the form
$\; H=\left(\begin{array}{cc}0&1\\1&0
\end{array} \right)\;$. The intersection form of the manifold
$S^{2}\times S^{2}$ realizes the form
$H$, and the K3 surface (quadric in ${\C}{\P}^{3}$) realizes
$ 2E_{8}\oplus 3H$. Donaldson had shown that if $k=1$, then $l\geq 3$
([{\bf D}]).
Clearly connected sums of K3 surface realizes
$\;2kE_{8}\oplus 3k H$. In general it is a conjecture that in (40) we must
necessarily have
$l\geq 3k$ (sometimes this is called $11/8$ conjecture).  Recently  by using
Seiberg-Witten theory M.Furuta has shown that

\vspace{.15in}
\noindent{\bf Theorem (Furuta) }: Let $X$ be a simply connected closed
smooth spin \\
$4$-manifold with the intersection form $q_{X}=2k E_{8}\oplus lH  $, then
$l\geq 2k+1$
\vspace{.15in}

Proof: We will only sketch the proof of  $\;l\geq 2k$. We pick
$L\longrightarrow X$ to be the trivial bundle (it is characteristic since $X$
is
spin). Notice that the spinor bundles
$$V^{\pm}=P\times _{\rho_{\pm}}{\C}^{2}\longrightarrow X$$
 $\rho_{\pm}: x\longmapsto q_{\pm} x $ , are  quoternionic vector bundles.
That is,
there is an action  $j: V^{\pm}\to V^{\pm}$  defined by $[p,x]\to [p,xj]$,
which is
clearly well defined. This action commutes with
$$\partiali: \Gamma(V^{+})\longrightarrow \Gamma(V^{-}) $$ Let
$A_{0}$ be the trivial connection, and write $\;\pm A=A_{0}\pm i\;a\in
\cal{A}(L)$

\begin{eqnarray*}
\partiali _{A}(\psi j)&=&\sum \rho (e^{k}) \;[\;\nabla _{k} +i\ a\;]\; (\psi
j)=
\sum \rho (e^{k}) \;[\;\nabla _{k}(\psi)j+ \psi j ia\;] \\ &=& \sum
\rho (e^{k}) \;[\;\nabla _{k}(\psi)j-  ia \psi j\;]=
\partiali_{-A}(\psi)j
\end{eqnarray*}
$\Z_{4}\;$ action $\;(A , \psi)\longmapsto (-A , \psi j)\;$  on
$\;\Omega^{1}(X)\times \Gamma(V^{+})\;$  preserves the compact set
 $${\cal M}_{0}=\tilde {\cal M}\cap \; ker( d^{*}) \oplus \Gamma (V^{+})$$

\vspace{.05in}

\noindent where  $\;\;\tilde {\cal M}=\{ (a, \psi)\in
\Omega^{1}(X)\oplus \Gamma(V^{+})\;|\;
\partiali_{A}(\psi)=0 \;,\;\;\;\rho(F_{A}^{+})=\sigma(\psi) \;\} $ \\ For
example
from the local description of $\sigma $ in (2) we can check
$$\sigma (\psi j)=\sigma (z+jw)j = \sigma(-\bar{w} + j\bar{z})=-\sigma(\psi)=
-\;F_{A}^{+}=F_{-A}^{+}$$ This $\Z_{4}$ in fact extends to an action of the
subgroup $G$ of
$SU(2)$ which is generated by $<S^{1}\;,\;j>$, where $S^{1}$ acts trivially on
$\Omega^{*}$ and by complex multiplication on $\Gamma(V^{+})$, and
$j$ acts by $-1$ on $\Omega^{*}$ and by quaternionic multiplication on
$\Gamma(V^{+})$ In particular we get a
$G$-equivariant map
$\;\varphi = L + \theta : {\cal V}\to {\cal W} \;$ where:

$$ \varphi  : {\cal V}=ker(d^{*})\oplus
\Gamma(V^{+})\longrightarrow  {\cal W}= \Omega_{+}^{2}\oplus
\Gamma(V^{-}) $$
$$L=\left( \begin{array}{cc} d^{+} & 0 \\ 0 & \partiali \end{array}
\right)\;\;\;
\mbox{and}\;\;\;
\theta (a,\psi)=(\;\sigma (\psi)\;,\; a \psi\;)$$ with
$\varphi^{-1}(0)= {\cal M}_{0}$ and $\varphi (v)=L(v) + \theta (v)$ with
$L$ linear
Fredholm and $\theta $ quadratic. We apply the ``usual" Kuranishi technique (cf
[{\bf L}]) to obtain a finite dimensional local model
$\;V\longmapsto W\;$ for $\varphi$.

We let  ${\cal V}=\oplus V_{\lambda }$ and ${\cal W}=\oplus W_{\lambda }$,
where
$V_{\lambda }$ and $W_{\lambda }$ be $\lambda $ eigenspaces of
$L^{*}L :V\to V$ and $LL^{*}:W\to W$
 repectively. Since $L L^{*}$ is a multiplication by $\lambda $ on
$V_{\lambda}$, for $\lambda > 0$ we have isomorphisms $L
:V_{\lambda}\stackrel{\cong
}{\longrightarrow} W_{\lambda }$. Now pick
$\Lambda >0$ and consider projections:

$$\oplus_{\lambda\leq
\Lambda}W_{\lambda}\;\;\stackrel{p_{\Lambda}}{\longleftarrow}
\;\; W \;\;
\stackrel{1-p_{\Lambda}}{\longrightarrow}\;\;
\oplus_{\Lambda > \Lambda} W_{\lambda}$$ Consider the local diffeomorphism
$f_{\Lambda}: V \stackrel{ }{\longrightarrow} V$ given by:
$$u=f_{\Lambda}(v)=v+L^{-1}(1-p_{\Lambda})\theta (v)
\;\;\Longleftrightarrow
\;\;L(u)=L(v)+(1-p_{\Lambda})\theta (v)$$ The condition
$\varphi(v)=0$ is equivalent to $\;p_{\Lambda}\;\varphi (v)=0 $ and
$ (1-p_{\Lambda})\;\varphi (v)=0 $, but
\begin{eqnarray*} (1-p_{\Lambda})\;\varphi (v)=0\;\Longleftrightarrow \;
(1-p_{\Lambda})\;L(v) + (1-p_{\Lambda})\;\theta(v)=0
\;\Longleftrightarrow\;\;\;\;\;\;
\;\;\;\;\;\;\;\;\;\;\; \\ (1-p_{\Lambda})\;L(v) + L(u)-L(v) =0
\; \Longleftrightarrow\;L(u)=p_{\Lambda}\;L(v) \;
\Longleftrightarrow\; u\in
\oplus_{\lambda\leq \Lambda}V_{\lambda}
\end{eqnarray*} Hence $\;\varphi
(v)=0\;\Longleftrightarrow\;p_{\Lambda}\;\varphi
(v)=0\; $ and
$\;u\in \oplus_{\lambda\leq \Lambda}V_{\lambda}\;$, let
$$ \varphi _{\Lambda}:V=\oplus_{\lambda\leq \Lambda}
V_{\lambda}\longrightarrow W=
\oplus_{\lambda\leq
\Lambda}W_{\lambda}\hspace{.15in}
\mbox{where}\hspace{.15in}
\varphi _{\Lambda}(u)=p_{\Lambda}\;\varphi\; f_{\Lambda}^{-1}(u)$$

Hence in the local diffeomorphism
$f_{\Lambda}: {\cal O}\stackrel{ \approx}{\longrightarrow} {\cal O}
\subset {\cal V }$ takes the piece of the compact set
$f_{\lambda}({\cal O}\cap {\cal M}_{0} )$ into the finite dimensional subspace
$V\subset {\cal V}$, where ${\cal O}$ is a neighborhood of $(0,0)$.
 As a side fact note that near $(0,0)$ we have
$${\cal M}(L) \approx{\cal M}_{0}(L)/S^{1} $$ We claim that for
$\lambda \gg 1 \;$,  the local diffeomorphism
$f_{\Lambda}: {\cal O}\stackrel{ \approx}{\longrightarrow} {\cal O}
\subset V $
 extends to a  ball $B_{R}$  of large radius $R$ containig the compact set
${\cal M}_{0}(L) $, i.e. we can make the zero set
$\varphi_{\Lambda}^{-1}(0)$ a compact set.

 We see this by applying the Banach contraction principle. For example for
a given
$u\in B_{R}$, showing that there is $v\in V$ such that
$f_{\Lambda}(v)=u$ is equivalent of showing that the map
$\;T_{u}(v)=u- L^{-1}(1-p_{\Lambda})\theta (v)\;$ has a fixed point. Since
$L^{-1}(1-p_{\Lambda})$ has eigenvalues
$1/\lambda$ on each
$W_{\lambda}$ in appropriate Sobolev norm we can write
$$|| T_{u}(v_{1})- T_{u}(v_{2}) ||\leq \frac{C}{\Lambda }||\theta(v_{1})-
\theta (v_{2})||
\leq \frac{C}{\Lambda }||v_{1}-v_{2}||$$

Vector subspaces $V_{\lambda}$ and $W_{\lambda}$ are either quaternionic or
real
depending on whether they are subspaces of
$\Gamma (V^{\pm})$ or
$\Omega^{*}(X)$. For a generic metric we can make  the cokernel of
$\partiali$ zero hence the dimension of the kernel (as a complex vector
space) is
$\;ind(\partiali)=-\sigma /8=2k\;$, and since $H^{1}(X)=0$  the dimension of
the
cokernel of $d^{+}$ (as a real vector space) is
$b^{+}=l$.
 Hence $\;\varphi _{\Lambda}\;$ gives  a $G$-equivariant map
$$\varphi: {\H}^{k+y}\oplus {\R}^{x}\longrightarrow {\H}^{y}\oplus
{\R}^{l+x}$$ with
compact zero set. From this Furuta shows that $l \geq 2l+1$. Here we give
an easier
argument of D.Freed which gives a slightly weaker result of $l\geq 2k$. Let
$E_{0}$ and $E_{1}$ be the complexifications of the domain an the range of
$\varphi
$; consider $E_{0}$ and $E_{1}$ as bundles over a point
$x_{0}$ with projections $\pi_{i}:E_{i}\to x_{0}$, and with
$0$-sections
$s_{i}:x_{0}\to E_{i}\;, i=0,1$. Recall $K_{G}(x_{0})=R(G)$, and we  have Bott
isomorphisms
$\beta (\rho)=\pi _{i}^{*}(\rho)\;\lambda_{E_{i}}\;$, for $\;i=0,1$ where
$\;\lambda_{E_{i}}\;$ are the Bott classes. By compactness  we get an
induced map
$\varphi ^{*}$:

$$\begin{array} {ccc} K_{G}(B(E_{1}),S(E_{1}))&\stackrel{\varphi ^{*}}
\longrightarrow & K_{G}(B(E_{0}),S(E_{0}))\\ &&\\
\approx \;\uparrow \beta & & \approx \;\uparrow \beta \\ R(G) && R(G)
\end{array}$$ Consider  $s_{i}^{*}(\lambda_{E_{i}})=
\sum(-1)^{k}\Lambda^{k}(E_{i})=\Lambda_{-1}(E_{i})\in R(G)$, then by some
$\rho$ we have
$$\Lambda_{-1}(E_{1})=s_{1}^{*}(\lambda_{E_{1}})=s_{0}^{*}
\varphi^{*}(\lambda_{E_{1}}) =s_{0}^{*} (\pi
_{0}^{*}(\rho)\;\lambda_{E_{0}})=\rho
\;\Lambda_{-1}(E_{0})$$ So in particular $tr_{j}(\Lambda_{-1}(E_{0}))$ divides
$tr_{j}(\Lambda_{-1}(E_{1}))$. By recalling $\;j:E_{i}\to E_{i}\;$
$$tr_{j}(\Lambda_{-1}(E_{i}))=det (I-j)\;\;\;\;\;\mbox{for}\;\;i=0,1$$ Since
$(z,w)j=(z+jw)j=-\bar{w}+j\bar{z}=(-\bar{w},\bar{z})\;$  $j$ acts on
${\H}\otimes{\C}$ by matrix
$$ A=\left(\begin{array}{cccc}0 &0 &1 &0\\ 0&0&0&-1\\-1&0&0&0\\0&1&0&0
\end{array}\right) $$ so $\;det (I-A)=4$, and $j$ acts on $\;{\R}\otimes
{\C}\;$ by
$j(x)=-x$
 so $det (I-(-I))=2\;$. Hence $\;4^{k+y}\;2^{x}\;$ divides
$\;4^{y}\;2^{l+x}\;$ which implies
$\;l\geq 2k\;$  \hspace{1.4in} $\Box$

\vspace{.15in}

There is another nice application of Seiberg-Witten invariants: It is an
old problem
whether the quotient of a simply connected smooth complex surface by an
antiholomorphic involution
$\sigma:\tilde{X}\to \tilde{X}\;$ (an involution which anticommutes with
the almost
complex homomorphism
$\sigma _{*}J=-J\sigma_{*}$) is a "standard" manifold (i.e. connected sums of
$S^{2}\times S^{2}$ and $\;\pm \C\P^{2}\;$). A common example of a
antiholomorphic involution is the complex conjugation on a complex projective
algebraic surface with real coefficients. It is known that the quotient of
$\C\P^{2}$ by complex conjugation is $S^{4}$ (Arnold, Massey, Kuiper); and
for every
$d$ there is a curve of degree $d$ in  $\C\P^{2}$ whose two fold
branched cover has a standard quotient ([{\bf A}]). This problem makes
sense only if
the antiholomorphic involution has a fixed point, otherwise the quotient
space has
fundamental group
$\Z_{2}$ and hence it can not be standard. By ``connected sum" theorem,
Seiberg-Witten invariants of ``standard" manifolds vanish, so it is natural
question
to ask whether
 Seiberg-Witten invariants of the quotients vanish. Shugang Wang has shown
that this
is the case for free antiholomorphic involutions.

\vspace{.15in}

\noindent {\bf Theorem} ({\bf S.Wang}) Let $\tilde{X}$ be a minimal Kahler
surface
of general type, and $\sigma: \tilde{X}\to \tilde{X}$ be a free antiholomorphic
involution, then the quotient
$X=\tilde{X}/\sigma$ has all Seiberg-Witten invariants zero

\vspace{.12in} Proof: Let $h$ be the Kahler metric on $\tilde{X}$, i.e.
$\omega(X,Y)=h(X,JY)$ is the Kahler form. Then
$\;\tilde{g}=h+\sigma^{*}h\;$ is an invariant metric on $\tilde {X}$ with the
Kahler form
$\;\tilde{\omega}=\omega - \sigma^{*}\omega \;$. Let $g$ be the
``push-down" metric
on $X$. Now we claim that all $SW_{L}(X)=0$ for all $L\to X$, in fact we
show that
there are no solutions to Seiberg-Witten equations for $X$: Otherwise if
$L\longrightarrow (X,g)$ is the characteristic line bundle supporting a
solution
$(A,\psi)$, then the pull-back pair
$(\tilde{A},\tilde{\psi})$ is a solution for the pull-back line bundle
$\tilde{L}\longrightarrow \tilde{X}$ with the pull-back $Spin_{c}$
structure, hence

$$0\leq  dim {\cal M_{\tilde{L}}}(\tilde{X})=\frac{1}{4}c_{1}^{2}(\tilde{L})-
\frac{1}{4} (3\sigma(\tilde{X})+2 \chi (\tilde{X}))$$ But
$\tilde{X}$ being a minimal Kahler suface of general type
$3\sigma(\tilde{X})+2 \chi (\tilde{X})=K^{2}_{\tilde{X}}>0\;$, hence
$\;c_{1}^{2}(\tilde{L})>0  $. This implies that
$(\tilde{A},\tilde{\psi})\;
$ must be an irreducible solution (i.e. $\psi \neq 0 $), otherwise
$F^{+}_{\tilde{A}}=0$ would imply
$\;c_{1}^{2}(\tilde{L})<0\;$. Now by (25) the nonzero solution
$\;\psi=
\alpha u_{0} + \beta$ must have either one of ${\;\alpha\;}$ or
$\;\beta\;$ is zero (so the other one is nonzero), and since
$\tilde{\omega}\wedge \tilde{\omega}$ is the volume element:
$$\tilde{\omega}.c_{1}(\tilde{L})=\frac{i}{2\pi}
\int \tilde{\omega}\wedge F^{+}_{\tilde{A}}=
\frac{i}{2\pi}\int \tilde{\omega}\wedge(\frac{|\beta |^{2}- |\alpha
|^{2}}{2})\;i\;\tilde{\omega}
\neq 0$$ But since  $\;\sigma^{*} (\tilde{\omega})=-
\tilde{\omega}\;$,
$\;\sigma^{*} c_{1}(\tilde{L})= c_{1}(\tilde{L})\;$, and  $\;\sigma
\;$ is an orientation preserving map we get a contradiction \\
$$\;\;\;\;\;\;\;\;\;\;\;\;\;\;\;
\tilde{\omega}.c_{1}(\tilde{L})=\sigma^{*}(\tilde{\omega}.c_{1}(\tilde{L}))=
-\tilde{\omega}.c_{1}(\tilde{L})
\;\;\;\;\;\;\;\;\;\;\;\;\;\;\;\;\;\;\;\;\;\;\;\Box$$

\end{document}